\documentclass[english]{article}
\usepackage[T1]{fontenc}
\usepackage[latin9]{inputenc}
\usepackage{amsmath}
\usepackage{amssymb}
\usepackage{stackrel}

\makeatletter
\@ifundefined{date}{}{\date{}}
\makeatother

\usepackage{babel}
\begin{document}
\title{A proof of Gödel's incompleteness theorems using Chaitin's incompleteness
theorem}
\author{David O. Zisselman}
\maketitle

\section{Abstract }

Gödel's first and second incompleteness theorems are corner stones
of modern mathematics. In this article we present a new proof of these
theorems for $ZFC$ and theories containing $ZFC$, using Chaitin's
incompleteness theorem and a very basic numbers extension.

As opposed to the usual proofs, these proofs don't use any fixed point
theorem and rely solely on sets structure. Unlike in the original
proof, the statements which can be shown to be unprovable by our technique
exceed by far one specific statement constructed from the axiom set.

Our goal is to draw attention to the technique of number extensions,
which we believe can be used to prove more theorems regrading the
provability and unprovability of different assertions regarding natural
numbers.

\section{Introduction}

Kolmogorov Complexity is the function $Kol\left(x\right)$ (where
$x\in\mathbb{N}$) which denotes the minimal $y$ such that $y$ is
a coding of a machine that calculates $x$. Given a consistent and
effective theory $T$, Chaitin's incompleteness theorem (cf: \ref{subsec:Chaitin's-incompleteness-theorem})
states that $T$ doesn't proof $Kol\left(x\right)\geq L$ for a large
enough $L$ (depending on $T$). Applying this to $T=ZFC$ and a fixed
$x$, one can construct a $ZFC$ model s.t $Kol\left(x\right)<L_{ZFC}$.
As it is impossible that 
\[
\forall x\in\mathbb{N}\,\,Kol\left(x\right)<L_{ZFC}
\]
 for that $L_{ZFC}$ we conclude that there must be a $x'\in\mathbb{N}$
for which $Kol\left(x\right)\geq L_{ZFC}$ holds true in $\mathbb{N}$
but can't be proven in $ZFC$ (and infinitely many such $x'$).

Moreover, if $ZFC$ could prove $con\left(ZFC\right)$, then in each
model of $ZFC$ one could continue this construction, i,e build a
series of models $M_{1},M_{2},\ldots$ of $ZFC$ and a series of natural
numbers $x_{1},x_{2},\ldots\in\mathbb{N}$ such that 
\[
\forall i\,\forall j\leq i\,\,M_{i}\models Kol\left(x_{j}\right)<L_{ZFC}
\]
 This construction is ``compressing'' by it's nature: For $i=L_{ZFC}+2$
we get that we have $L_{ZFC}+2$ \textbf{different }numbers $y_{1},y_{2},\ldots,y_{L+2}$
which represent machines which calculates $x_{1},x_{2},\ldots,x_{L+2}$.
As every $y_{i}\,|\,\,1\leq i\leq L+2$ held $y_{i}<L_{ZFC}$ we get
that the model $M_{L+2}$ has $L_{ZFC}+2$ \textbf{different }natural
numbers all smaller than $L_{ZFC}$, which is a contradiction. 

\section{Preliminaries}

Since the proofs use a few known theorems and basic definitions, I
shall quote the theorems which I'll use later. 

\subsection{Definitions and notations \label{subsec:Definitions-and-notations}}

For the following definitions we will assume we have one universal
Turing machine (i.e a coding method of Turing machines), by which
we measure the length of other Turing machines 
\begin{itemize}
\item Kolmogorov Complexity (single value) - The Kolmogorov Complexity of
a natural number $x$ will be denoted by $Kol\left(x\right)$.
\end{itemize}

\paragraph{Notation:}

A \textbf{Turing machine (TM)} will refer to a RAM computational machine. 

\paragraph{Definition:}

A \textbf{structure }for the language of set theory is a set $A$
along with a two place relation $\left(\in\right)\subset A^{2}$ . 

I.e a set of elements $A$ along with a subset $R\subset A^{2}$ which
we will denote by $a\in b\Leftrightarrow aRb$ 

\paragraph{Definition ($ZF$):}

the set of axioms on the language of set theory called $ZF$ (Zermelo--Fraenkel
axiomatic set theory) is the following set:
\begin{enumerate}
\item Axiom of extensionality $\forall x\forall y\left(\left(\forall z\left(z\in x\leftrightarrow z\in y\right)\rightarrow x=y\right)\right)$
\item Axiom of regularity $\forall x\left(\exists a\left(a\in x\right)\rightarrow\exists y\left(y\in x\wedge\neg\exists z\left(z\in y\wedge z\in x\right)\right)\right)$
\item Axiom schema of specification $\forall z\forall w_{1}\forall w_{2}...\forall w_{n}\exists y\forall x\left(x\in y\leftrightarrow\left(x\in z\wedge\phi\left(w_{1},...,w_{n},x\right)\right)\right)$
\item Axiom of pairing $\forall x\forall y\exists z\left(x\in z\wedge y\in z\right)$
\item Axiom of union $\forall F\exists A\forall Y\forall x\left(\left(x\in Y\wedge Y\in F\right)\rightarrow x\in A\right)$
\item Axiom schema of replacement 
\[
\!\!\!\!\!\!\!\!\!\!\!\!\!\!\!\!\!\!\!\!\!\!\!\!\!\!\!\!\!\!\!\!\!\!\!\forall A\forall w_{1}\forall w_{2}...\forall w_{n}\left[\forall x\left(x\in A\rightarrow\exists!y\phi\left(w_{1},...,w_{n},x,y\right)\right)\rightarrow\exists B\forall x\left(x\in A\rightarrow\exists y\left(y\in B\wedge\phi\left(w_{1},...,w_{n},x,y\right)\right)\right)\right]
\]
\item Axiom of infinity $\exists X\left(\emptyset\in X\wedge\forall y\left(y\in X\rightarrow S\left(y\right)\in X\right)\right)$
\item Axiom of power set $\forall x\exists y\forall z\left(z\subseteq x\rightarrow x\in y\right)$
\end{enumerate}
\begin{itemize}
\item Definition: Axiom of choice (AC) is the following statement:
\[
\forall X\left[\emptyset\not\in X\Rightarrow\exists f:X\rightarrow\cup X\,\,\,\forall a\in X\left(f\left(a\right)\in a\right)\right]
\]
\item Definition: the set of axioms $ZF$ with AC is called $ZFC$ axioms
set.
\item Definition: the set of axioms $ZF$ without regularity is denoted
$ZF^{-}$ 
\item Definition: a definable language is a set of the form $\left\{ x\in\mathbb{N}\vert\phi\left(x\right)\right\} $
where $\phi$ is a one free variable formula in the language of $ZF$. 
\end{itemize}
Unless mentioned otherwise, a language will refer only to a definable
language.

\paragraph{Definition:}

A model of $ZF$ or $ZFC$ or $ZF^{-}$ etc is a structure for the
language of set theory which holds the appropriate axioms of $ZF$
or $ZFC$ or $ZF^{-}$ respectively. 

\paragraph{Definition and notation of~ $\omega$}

Within a model $M$ of $ZF$ the standard way to construct the ``natural
numbers'' within the model is to define 0 to be $\emptyset$ and
$S\left(x\right)=x\cup\left\{ x\right\} $ (successor operation).
So $1=\left\{ \emptyset\right\} $ and $2=\left\{ \left\{ \emptyset\right\} ,\emptyset\right\} $
and so on. The minimal set which contains $\emptyset$ and is closed
under $S$ operation is called $\omega$ and the existence of it is
guaranteed by the axiom of infinity. However, and even though we think
about $\omega$ as $\mathbb{N}$ they may be very different objects
in some models. 
\begin{itemize}
\item The notation of $\omega\left(M\right)$ will be a more suitable notation
instead of the imprecise notation $\mathbb{N}$.
\item Given a model $M$ of $ZF$ denote $\omega\left(M\right)$ to be the
$\omega$ of the model $M$.
\item Given an $n\in\omega\left(M\right)$ denote 
\[
\left[n\right]\triangleq\left\{ i\in\omega\left(M\right)|i\leq n\right\} =\left\{ 0,1,2,\ldots,n\right\} 
\]
\item Unless specifically mentioned otherwise $\mathbb{N}$ will denote
the set of natural numbers $\mathbb{N}=\left\{ 0,1,2,3,4,...\right\} $
which includes zero.\\
As, strictly speaking, $\mathbb{N}$ can't be defined I will not use
the notation $\mathbb{N}$ except for intuitions. 
\end{itemize}

\paragraph{Definition (set model):}

Given a ZF model $V_{1},\in_{V_{1}}$, the model $V_{2},\in_{V_{2}}$
of ZF is said to be \textbf{inside} $V_{1}$ (or a set of $V_{1}$)
if the following hold:
\begin{itemize}
\item $V_{2}$ is a set within $V_{1}$ 
\item all sets of $V_{2}$ are also sets in $V_{1}$.
\item $\in_{V_{2}}$ is a set of 2-tuples (i.e as a subset of $\left(V_{2}\right)^{2}$)
is a set of $V_{1}$
\end{itemize}

\paragraph{Lemma:}

Let $V_{1}$ and $V_{2}$ be $ZF$ models. If $V_{2}$ is a set of
$V_{1}$ and if $a$ is a set of $V_{2}$ and $b\in_{V_{2}}a$ , then
$b$ is a set of $V_{2}$ and $b$ is a set of $V_{1}$.

\paragraph{Proof:}

Under the above condition~ $b\in_{V_{2}}a$. As $\in_{V_{2}}$ is
a subset of $\left(V_{2}\right)^{2}$, this means $b$ is a set of
$V_{2}$. As $V_{2}$ is inside $V_{1}$, we get that $b$ is also
a set of $V_{1}$. $\blacksquare$

\subsubsection{Non - compression lemma \label{subsec:non---compression}}

\paragraph{Lemma:}

Given a model $M$ of $ZF$ and given $x\in\omega\left(M\right)$,
there exist~ $y\in\omega\left(M\right)$~ such that 
\[
Kol\left(y\right)>x
\]

\paragraph{Proof:}

The reader may want to think first about the case~ $\omega\left(M\right)=\mathbb{N}$
first.

The function that maps $y$ to $z$, the smallest code of TM that
prints $y$ , is a definable function in $ZF$. Assume by negation
there doesn't exist $y\in\omega\left(M\right)$ such that $Kol\left(y\right)>x$
. In this case 
\[
M\vDash\forall y\,\,Kol\left(y\right)\leq x
\]
 In such a case, there is a mapping from $\left[x\right]=\left\{ 0,1,,\ldots,x\right\} \subset\omega\left(M\right)$
to $\omega\left(M\right)$ which is both injective and surjective
(as for every $y$ exist $z\in\left[x\right]$ which codes a TM which
calculates $y$) and is definable within $M$. Hance, within $M$,
the set $\omega\left(M\right)$ is of cardinally $\leq x$ and $x<\aleph_{0}$
, which is a contradiction. $\blacksquare$

\paragraph{Notation:}

given a model $M$ of $ZF$ and given $x\in\omega\left(M\right)$
define $BB_{M}\left(x\right)$ to be the smallest number $y\in\omega\left(M\right)$
s.t $M\vDash Kol\left(y\right)>x$. 

\subsubsection{Chaitin's incompleteness theorem (1971) \label{subsec:Chaitin's-incompleteness-theorem} }

\paragraph{Theorem:}

Let $V$ be a model of $ZF$. Let $T$ be an effective consistent
set of axioms. Then, there exist $L\in\mathbb{\omega}\left(V\right)$
(which depends on the set of axioms) such that for every $x\in\mathbb{\omega}\left(V\right)$
the statement $L\leq Kol\left(x\right)$ can't be proven from the
set $T$ .

For completeness of the article I'll add proof to this theorem:

The reader may want to consider the case of $\omega\left(V\right)=\mathbb{N}$
at first.

\paragraph{Proof: }

Within $V$, any proof of claim $\phi$ from $T$ is a number within
$\mathbb{\omega}\left(V\right)$. Denote $w\in\omega\left(V\right)$
to be the first proof of a claim in $\left\{ L\leq Kol\left(x\right)|x\in\omega\left(V\right)\right\} $
and let $x'\in\mathbb{\omega}\left(V\right)$ be the number s.t $T\underset{w}{\vdash}L\leq Kol\left(x'\right)$
(the proof $w$ proves that $L\leq Kol\left(x'\right)$). Create a
Turing machine $M$ that goes over $y\in\mathbb{\omega}\left(V\right)$
(in the regular order) and checks if $y$ is a proof of the statement
$L\leq Kol\left(x\right)$ and halts if it is and prints $x$. The
size of the TM $M$ (as TM could be coded as natural numbers too)
is $\log\left(L\right)+C$ for a fixed $C\in\mathbb{\omega}\left(V\right)$.
As, the decimal representation of $L$ is $\log\left(L\right)$ and
$C$ is the extra size required to represent the theory $T$ and the
operation of $M$. 

The TM $M$ is a TM that prints $x'$ and its size is $\log\left(L\right)+C$
and as $L\leq Kol\left(x'\right)$ it follows that $L\leq\log\left(L\right)+C$
(as we've just shown a TM of size $\log\left(L\right)+C$). the inequality
$L\leq\log\left(L\right)+C$ can't hold for a sufficiently large $L$. 

\subparagraph{notation:}

given an effective consistent set of axioms $T$ define $\mathcal{L}\left(T\right)$
to be the minimal number $L$ s.t $L>\log\left(L\right)+C$. 

\paragraph{Remarks:}
\begin{itemize}
\item Notice that $\mathcal{L}\left(T\right)\leq2C$.
\item The notation $\mathcal{L}\left(T\right)$ doesn't mention $V$ because
the value doesn't depend on $V$ to a large extent (and it will later
be stated and proved).
\end{itemize}

\section{Extension}

On this part I'll describe the basic principles of numbers extensions.

\subsubsection{model extension definition \label{subsec:model-extension-definition}}

\paragraph{Definition: }

Given two models $V_{1},V_{2}$ of $ZF$ we say that $V_{2}$ number
extends $V_{1}$ using mapping $f:\omega\left(V_{1}\right)\rightarrow\omega\left(V_{2}\right)$
if the arithmetic operations of $V_{1},V_{2}$ i,e 
\[
+_{1},+_{2},\cdot_{1},\cdot_{2}\wedge_{1},\wedge_{2}
\]
accordingly are respected by $f$. Meaning: 
\begin{enumerate}
\item $\forall a,b\in\omega\left(V_{1}\right)\left(f\left(a+_{1}b\right)=f\left(a\right)+_{2}f\left(b\right)\right)$
\item $\forall a,b\in\omega\left(V_{1}\right)\left(f\left(a\cdot_{1}b\right)=f\left(a\right)\cdot_{2}f\left(b\right)\right)$
\item $\forall a,b\in\omega\left(V_{1}\right)\left(f\left(a^{b}\right)=f\left(a\right)^{f\left(b\right)}\right)$
\item $f\left(0_{V_{1}}\right)=0_{V_{2}}$and $f\left(1_{V_{1}}\right)=1_{V_{2}}$
\end{enumerate}
I'll denote it by $V_{1}\overset{f}{\longrightarrow}V_{2}$. 

\subparagraph{Remark:}
\begin{itemize}
\item Please note that there is no assumption that inductive arguments on
elements of $\omega\left(V_{1}\right)$ are ``transferred'' (in
some way or another) into $\omega\left(V_{2}\right)\cap f\left(\omega\left(V_{1}\right)\right)$.
\item The function $f$ itself isn't assumed to be a set of either $V_{1}$
or $V_{2}$.
\end{itemize}

\subsubsection{Model initial segment extension definition}

\paragraph{Definition:}

If $V_{1}\overset{f}{\longrightarrow}V_{2}$. and the extension holds
\[
\forall b\in\omega\left(V_{1}\right)\forall a\in\omega\left(V_{2}\right)\exists c\in\omega\left(V_{1}\right)\left(a<f\left(b\right)\rightarrow a=f\left(c\right)\right)
\]
 the extension is said to be \textbf{initial segment extension} (i.s
extension) and will be denoted by $V_{1}\stackrel[I]{f}{\rightarrow}V_{2}$
.

\paragraph*{Please note the following:}
\begin{enumerate}
\item The meaning of i.s extension is that the model $V_{2}$ has ``more
numbers'' than $V_{1}$ but the numbers $V_{1}$ ``behave the same''
on $V_{2}$ ``lower parts''.
\item The function $f$ is usually undefinable from the models themselves.
\item the fact that $V_{1}\overset{f}{\rightarrow}V_{2}$ for two models
doesn't imply that the two models are elementary equivalent. Nor does
it mean they agree on $\omega$ attributed formulas.
\end{enumerate}

\subsection{Basic properties of extension\label{subsec:basic-properties-of-extension}}

\paragraph{Lemma:}

Given two models of $ZF$ , $V_{1},V_{2}$ s.t 
\begin{itemize}
\item $V_{2}$ is a set within $V_{1}$ 
\item $\in_{2}$ is a set (of pairs) within $V_{1}$
\end{itemize}
then a mapping $f\in V_{1}$, $f:\omega\left(V_{1}\right)\rightarrow\omega\left(V_{2}\right)$
exist s.t $V_{1}\underset{I}{\overset{f}{\longrightarrow}}V_{2}$.

\paragraph{Proof:}

The construction of $f$ is by induction:
\begin{itemize}
\item The base case: define $f\left(0_{v_{1}}\right)=0_{V_{2}}$.
\item For $a=b+1$ where $a,b\in_{V_{1}}\omega\left(V_{1}\right)$ define
$f\left(a\right)=f\left(b+_{1}1\right)=f\left(b\right)+_{V_{2}}1$
\item As every natural number is either $0$ or $b+1$ for another natural
number $b$, $f$ is defined.
\end{itemize}
Recall that the definitions of $+,\times,\wedge$ are inductive as
well. For $a,c\in_{V_{1}}\omega\left(V_{1}\right)$ 
\begin{itemize}
\item If $c=0_{V_{1}}$ then by definition:
\begin{itemize}
\item $a+c=a+0=a$
\item $a\cdot c=a\cdot0=0$
\item $a^{c}=a^{0}=1$ and if $a\neq0$ then $c^{a}=0^{a}=0$
\end{itemize}
\item Therefore if $c=0$, $f$ holds the equalities in \ref{subsec:model-extension-definition}:
\begin{itemize}
\item $f\left(a+_{1}c\right)=f\left(a+_{1}0\right)=f\left(a\right)=f\left(a\right)+_{2}0_{V_{2}}$.
\item $f\left(a\cdot_{1}c\right)=f\left(a\cdot_{1}0_{V_{!}}\right)=f\left(0_{V_{1}}\right)=0_{V_{2}}=f\left(a\right)\cdot_{2}0_{V_{2}}=f\left(a\right)\cdot_{2}f\left(c\right)$
\item $f\left(a^{c}\right)=f\left(a^{0}\right)=f\left(1_{V_{1}}\right)=1_{V_{2}}=f\left(a\right)^{0_{V_{2}}}=f\left(a\right)^{f\left(c\right)}$
\item If $a\neq0$ then $f\left(c^{a}\right)=f\left(0^{a}\right)=f\left(0_{V_{1}}\right)=0_{V_{2}}=f\left(c\right)^{f\left(a\right)}$
\end{itemize}
\item For $c\neq0$ exist $b\in_{V_{1}}\omega\left(V_{1}\right)$ s.t $c=b+1$
and thus by definition:
\begin{itemize}
\item $a+c=a+\left(b+1\right)=\left(a+1\right)+b$
\item $a\cdot c=a\cdot\left(b+1\right)=\left(a\cdot b\right)+a$
\item $a^{c}=a^{b+1}=a\cdot\left(a^{b}\right)$
\end{itemize}
\item Therefore if $c\neq0$, $f$ holds the equalities in \ref{subsec:model-extension-definition}:
\begin{itemize}
\item $f\left(a+_{1}c\right)=f\left(\left(a+_{1}1\right)+_{1}b\right)=\left(f\left(a\right)+_{2}1_{V_{2}}\right)+_{2}f\left(b\right)$.
\item $f\left(a\cdot_{1}c\right)=f\left(\left(a\cdot_{1}b\right)+_{1}a\right)=f\left(a\cdot b\right)+_{2}f\left(a\right)=f\left(a\right)\cdot_{2}f\left(b\right)+_{2}f\left(a\right)=f\left(a\right)\left(f\left(b\right)+_{2}1_{V_{2}}\right)=f\left(a\right)f\left(b+_{1}1\right)=f\left(a\right)\cdot_{2}f\left(c\right)$
\item $f\left(a^{c}\right)=f\left(a\cdot_{1}\left(a^{b}\right)\right)=f\left(a\right)\cdot_{2}f\left(a^{b}\right)=f\left(a\right)\cdot_{2}f\left(a\right)^{f\left(b\right)}=f\left(a\right)^{f\left(b\right)+1_{V_{2}}}=f\left(a\right)^{f\left(c\right)}$
. 
\end{itemize}
\end{itemize}
So far we've seen that $V_{1}\overset{f}{\rightarrow}V_{2}$. 

This $f$ also holds the extension property as well: 

As $V_{2}$ is a set of $V_{1}$, in $V_{1}$ one can define the following
set 
\[
S=\left\{ b\in\omega\left(V_{1}\right)\,|\,\exists a\in\omega\left(V_{2}\right)\forall c\in\omega\left(V_{1}\right)\left(\left(a<_{2}f\left(b\right)\right)\wedge\left(a\neq f\left(c\right)\right)\right)\right\} 
\]
if $S$ isn't empty and as $S$ is a set of natural numbers in $\omega\left(V_{1}\right)$
it has a minimum. Denote 
\[
s=\min S
\]
 notice that $s\neq0$ as $f\left(0_{V_{1}}\right)=0_{V_{2}}$ and
$\neg\exists a\in\omega\left(V_{2}\right)\,\,\left(a<0_{V_{2}}=f\left(0_{V_{1}}\right)\right)$.
So $s>0$ and as such $s=m+1$ for some $m\in_{V_{1}}\omega\left(V_{1}\right)$.
Therefore $f\left(s\right)=f\left(m\right)+_{2}1_{V_{2}}$. As $s$
was minimal it holds that $m\not\in_{V_{1}}S$. and thus for $m$
\begin{equation}
\forall a\in\omega\left(V_{2}\right)\exists c\in\omega\left(V_{1}\right)\left(a<f\left(m\right)\rightarrow a=f\left(c\right)\right)
\end{equation}
 as $s\in S$ we know 
\[
\exists a\in\omega\left(V_{2}\right)\forall c\in\omega\left(V_{1}\right)\left(\left(a<_{2}f\left(s\right)\right)\wedge\left(a\neq f\left(c\right)\right)\right)
\]
let $\exists a\in\omega\left(V_{2}\right)$ be constant and receive
\begin{equation}
\forall c\in\omega\left(V_{1}\right)\left(\left(a<_{2}f\left(s\right)\right)\wedge\left(a\neq f\left(c\right)\right)\right)
\end{equation}

if $a<_{2}f\left(m\right)$ then by (1) we know that $\exists c\in\omega\left(V_{1}\right)$
that violates (2). Otherwise if $a=_{2}f\left(m\right)$ condition
(2) is violated with $c=m$ . Lastly if $a>_{2}f\left(m\right)$ it
holds that $a\geq_{2}f\left(m\right)+1_{V_{2}}=f\left(s\right)$ which
is a contradiction to $a<_{2}f\left(s\right)$ in (2). As we've received
a contradiction it must be the case that $S$ above is empty. Thus
\[
\forall b\in\omega\left(V_{1}\right)\forall a\in\omega\left(V_{2}\right)\exists c\in\omega\left(V_{1}\right)\left(a<f\left(b\right)\rightarrow a=f\left(c\right)\right)
\]

$\blacksquare$

\paragraph{Remarks:}

Please note that the condition $V_{2}$ is a set of $V_{1}$ give
us $V_{1}\underset{I}{\overset{f}{\longrightarrow}}V_{2}$ but in
fact it's a much much stronger assertion. As in this case $f$ is
also in $V_{1}$.

Specifically, one could use this fact to define induction on $\left[a\right]$
for $a\in\omega\left(V_{2}\right)\cap Im\left(f\right)$ and pulls
back the argument by $f^{-1}$ to an induction on $\omega\left(V_{1}\right)$.
And, use this fact to use inductive arguments on $Im\left(f\right)$
as seen above.

\subsubsection{Absoluteness of Turing machines in sub-models\label{subsec:absoluteness-of-Turing}}

\paragraph{Theorem:}

Let $V_{1},V_{2}$ be two models of $ZF$ s.t $V_{2}$ is a set of
$V_{1}$ ,$V_{1}\underset{I}{\overset{f}{\longrightarrow}}V_{2}$
and let $T_{1}\in\omega\left(V_{1}\right)$ represent a coding of
a Turing machine. Let $z\in\omega\left(V_{1}\right)$ be any number. 

Denote $R\left(T_{1},z\right)\in\omega\left(V_{1}\right)$ to be the
coded state of the machine $T_{1}$ after $z$ steps (in $V_{1}$)
and $R\left(f\left(T_{1}\right),f\left(z\right)\right)$ be the coded
state of the machine $f\left(T_{1}\right)$ after $f\left(z\right)$
steps (in $V_{1}$).

Then 
\[
f\left(R\left(T_{1},z\right)\right)=R\left(f\left(T_{1}\right),f\left(z\right)\right)
\]

\paragraph{Proof (sketch):}

Denote $J$ to be the operation of running the machine in a specific
status one more step, i,e the function that 
\[
J\left(R\left(T,z\right)\right)=R\left(T,z+1\right)
\]
 for any TM $T\in\omega\left(V_{1}\right)$ and any number $z\in\omega\left(V_{1}\right)$.
The function $J$ can be expressed using arithmetic operations. The
proof is done by induction on $\omega\left(V_{1}\right)$:

The base case is $R\left(T_{1},0_{V_{`}}\right)$ as $T_{1}$ is a
machine that is coded by a text to integer coding , it holds that
$f\left(T_{1}\right)$ is a text to integer coding of the same machine
in $\omega\left(V_{2}\right)$. Thus we get 
\[
f\left(R\left(T_{1},0_{V_{1}}\right)\right)=R\left(f\left(T_{1}\right),f\left(0_{V_{1}}\right)\right)=R\left(f\left(T_{1}\right),0_{V_{2}}\right)
\]
 for $z+1$ we recall that $J$ is an arithmetic function and thus
$f\left(J\left(x\right)\right)=J\left(f\left(x\right)\right)$ for
any $x\in\omega\left(V_{1}\right)$ and thus 
\begin{align*}
f\left(R\left(T,z+1\right)\right) & =f\left(J\left(R\left(T_{1},z\right)\right)\right)\\
=f\left(J\left(R\left(T_{1},z\right)\right)\right) & =J\left(f\left(R\left(T_{1},z\right)\right)\right)\\
=J\left(R\left(f\left(T_{1}\right),f\left(z\right)\right)\right) & =R\left(f\left(T_{1}\right),f\left(z\right)+1\right)\\
 & =R\left(f\left(T_{1}\right),f\left(z+1\right)\right)
\end{align*}
Now for the induction, denote the set 
\[
S=\left\{ z\in\omega\left(V_{1}\right)\,|\,f\left(R\left(T_{1},z\right)\right)\neq R\left(f\left(T_{1}\right),f\left(z\right)\right)\right\} 
\]
 the set $S$ is a subset of natural numbers definable in $V_{1}$
(as $V_{2}$ is a set of $V_{1}$). If $S$ isn't empty, it must have
a minimum let $z'$ be that minimum. $z'\neq0_{V_{1}}$ as we've shown
that 
\[
f\left(R\left(T_{1},0_{V_{1}}\right)\right)=R\left(f\left(T_{1}\right),0_{V_{2}}\right)
\]
 if $z'=z''+1$ then $z''\not\in S$ and thus 
\[
f\left(R\left(T_{1},z''\right)\right)=R\left(f\left(T_{1}\right),f\left(z''\right)\right)
\]
 and thus we get that $z'=z''+1\not\in S$ as 
\[
f\left(R\left(T_{1},z''+1\right)\right)=R\left(f\left(T_{1}\right),f\left(z''+1\right)\right)
\]
 Therefore, $S$ above must be empty and thus 
\[
f\left(R\left(T_{1},z\right)\right)=R\left(f\left(T_{1}\right),f\left(z\right)\right)
\]
for every $z\in\omega\left(V_{1}\right)$.$\blacksquare$

\paragraph{Remark:}
\begin{enumerate}
\item Please note that $V_{2}$ may may more numbers which are not in $Im\left(f\right)$.
In this case, it may be that for a TM $T_{1}$, $T_{1}$ will not
halt in $V_{1}$ but will halt in $V_{2}$ . How every for any step
in $Im\left(f\right)$ the running in $V_{1}$ and $V_{2}$ will agree.
\item The ``sketch'' part of this proof is the fact that the coding of
$R$ and the $J$ operation wasn't fully defined. As I trust the reader
is familiar with such constructions, I don't see added value in elaborating. 
\item Please note the importance of the assumption: $V_{2}$ is a set of
$V_{1}$. The fact $V_{1}\underset{I}{\overset{f}{\longrightarrow}}V_{2}$
alone isn't enough for this proof (as we need to use induction on
$\omega\left(V_{1}\right)$). 
\end{enumerate}

\subsubsection{Absoluteness of $\mathcal{L}\left(T\right)$ \label{subsec:absoluteness-of-L(T)}}

\paragraph{Theorem:}

Let $V_{1},V_{2}$ be two models of $ZF$ s.t $V_{2}$ is a set of
$V_{1}$, $V_{1}\underset{I}{\overset{f}{\longrightarrow}}V_{2}$
and let $T_{1}\in\omega\left(V_{1}\right)$ be a coding of a TM the
recognizes $T$ , an effective consistent set of axioms (we assume
here that $T$ is consistent according to both $\omega\left(V_{1}\right)$
and $\omega\left(V_{2}\right)$) .Let $\mathcal{L}_{1}\left(T_{1}\right),\mathcal{L}_{2}\left(f\left(T_{1}\right)\right)$
be $\mathcal{L}\left(T\right)$ computed within $V_{1},V_{2}$ respectively. 

Then $f\left(\mathcal{L}_{1}\left(T_{1}\right)\right)=\mathcal{L}_{2}\left(f\left(T_{1}\right)\right)$.

\paragraph{Intuition:}

Recall the definition of $\mathcal{L}\left(T\right)$ given in \ref{subsec:Chaitin's-incompleteness-theorem}:
$\mathcal{L}\left(T\right)$ to be the minimal number s.t $L>\log\left(L\right)+C$.
Therefore, as long as $C$ is interpreted the same in both models,
$\mathcal{L}\left(T\right)$ will also be the same. 

Moreover, $C$ contains:
\begin{enumerate}
\item Representation of a TM which computes $T$ (which is assumed to be
the same).
\item A representation of the machine $M$ which goes over the proofs from
$T$ and finds the first proof of the claim in $\left\{ L\leq Kol\left(x\right)|x\in\omega\left(V\right)\right\} $
(for a given fixed $L$).
\end{enumerate}
As both of these are absolute in sub-models (as seen in \ref{subsec:absoluteness-of-Turing})
$C$ must be interpreted the same.

\paragraph{Proof (sketch):}

The reader may want to think of the case $\omega\left(V_{1}\right)=\mathbb{N}$
at first. The idea of this proof is that all of \ref{subsec:Chaitin's-incompleteness-theorem}
construction could be done in a bounded set of $\omega\left(V_{1}\right)$
and $\omega\left(V_{2}\right)$ behaves ``the same'' within on such
subsets (lower parts). Here are some details:

Recall the proof of \ref{subsec:Chaitin's-incompleteness-theorem},
the proof creates a Turing machine $M$ which goes over over all elements
of $\omega\left(V_{1}\right)$ until it finds the first element $w$
which proves a claim in the set $\left\{ L\leq Kol\left(x\right)|x\in\omega\left(V\right)\right\} $.
\begin{itemize}
\item Proofs are a list of statements in first order logic, each statement
can be either an axiom or derived from the previous statements.
\item Such a proof can be coded into an integer be text to integer coding.
Such a coding require only the operations of addition, multiplication,
exponentiation.
\item The proof could be verified by a Turing machine, which runs the same
(by \ref{subsec:absoluteness-of-Turing}) in both $V_{1}$ and $V_{2}$.
\item Therefore the proofs interprets the same in both $V_{1}$ and $V_{2}$.
\begin{itemize}
\item i.e if $T\underset{p}{\vDash}\phi$ in $V_{1}$ then $f\left(T\right)\underset{f\left(p\right)}{\vDash}f\left(\phi\right)$
in $V_{2}$. where $T,p,\phi$ are coding of a theory, proof and a
statement respectively.
\end{itemize}
\item Let $w\in\omega\left(V_{1}\right)$ be the first proof of a statement
in $\left\{ L\leq Kol\left(x\right)|x\in\omega\left(V_{1}\right)\right\} $
in $V_{1}$. The proof $f\left(w\right)$ is also a proof of a statement
in $\left\{ L\leq Kol\left(x\right)|x\in\omega\left(V_{2}\right)\right\} $
in $V_{2}$ and due to the one to one correspondence of $f$, $f\left(w\right)$
must be the first such proof in $\omega\left(V_{2}\right)$ as well.
\item Therefore the operation of machine $M$ will work ``the same'' and
return $x$ in $V_{1}$ and $f\left(x\right)$ in $V_{2}$.
\item for $L_{1}\in\omega\left(V_{1}\right)$,$L_{2}\in\omega\left(V_{2}\right)$
we've created two Turing machines $M_{1}$ in $V_{1}$ and $f\left(M_{1}\right)$
in $V_{2}$ that computes $x',f\left(x'\right)$ in $V_{1},V_{2}$
respectively. Therefore in both $V_{1},V_{2}$ the two inequalities
must hold:
\[
L_{1}\leq\log\left(L_{1}\right)+C
\]
\[
L_{2}\leq\log\left(L_{2}\right)+f\left(C\right)
\]
 and this the minimum number that violates them must be ``the same''
i.e 
\[
f\left(\mathcal{L}_{1}\left(T_{1}\right)\right)=\mathcal{L}_{2}\left(f\left(T_{1}\right)\right)
\]
\end{itemize}

\section{Gödel's incompleteness theorems}

\subsection{First theorem}

\paragraph{Theorem (Gödel's first incompleteness):}

Let $T$ be an effective consistent set of axioms containing $ZF$
then there is an $\omega$ attributed claim $\phi$ which is independent
of $T$ (can't be proven not disproven)

\paragraph{Proof:}

As $T$ is consistent, let $M$ be a model of $T$. Denote within
$M$ 
\[
x=BB_{M}\left(\mathcal{L}\left(T\right)\right)
\]
and denote the claim $\phi$ to be 
\[
\mathcal{L}\left(T\right)<Kol\left(x\right)
\]
By Chaitin's theorem \ref{subsec:Chaitin's-incompleteness-theorem}
the claim 
\[
\mathcal{L}\left(T\right)<Kol\left(x\right)
\]
 can't be proven from $T$. However, in $M$, by the definition $BB_{M}$
and the non - compression lemma \ref{subsec:non---compression} it
holds that 
\[
M\models\mathcal{L}\left(T\right)<Kol\left(x\right)
\]
 so the claim $\phi$ can't be disproven either as $T$ has a model
that holds it to be true. $\blacksquare$

\subsection{Second theorem}

\paragraph{Theorem (Gödel's second incompleteness):}

Let $T$ be an effective consistent set of axioms containing $ZF$
then the claim $con\left(T\right)$ (i.e there is no proof of contradiction
from $T$) can't be proven from $T$.

\paragraph{Proof:}

Assume (by negation) that the claim $con\left(T\right)$ could be
proven from $T$. As $T$ is consistent, let $M_{0}$ be a model of
$T$. Construct a series of models $\left(M_{i}\right)_{i\in\omega\left(M_{0}\right)}$
in the following way: Given $M_{i}$ denote 
\[
x_{i}=BB_{M_{i}}\left(\mathcal{L}\left(T\right)\right)
\]
 within $M_{i}$. In $M_{i}$, by the definition $BB_{M_{i}}$ and
the non - compression lemma \ref{subsec:non---compression} it holds
that 
\[
M_{i}\models\mathcal{L}\left(T\right)<Kol\left(x_{i}\right)
\]
As $M_{i}$ proves that $T$ is consistent (here we use our assumption
that $con\left(T\right)$ could be proven from $T$) by Chaitin's
theorem \ref{subsec:Chaitin's-incompleteness-theorem} the claim 
\[
\mathcal{L}\left(T\right)<Kol\left(x_{i}\right)
\]
 can't be proven from $T$. So there exist a model $M_{i+1}$ (within
$M_{i}$) s.t 
\[
M_{i+1}\models\mathcal{L}\left(T\right)\geq Kol\left(x_{i}\right)
\]

denote $y_{i}\in_{M_{i+1}}\omega\left(M_{i+1}\right)$ to be the code
of the TM which calculates the value of $x_{i}$ in $M_{i+1}$.

As $\mathcal{L}\left(T\right)$ is absolute (by \ref{subsec:absoluteness-of-L(T)})
we know that in all models $M_{i}$ it is the same. As $M_{i+1}$
is a set of $M_{i}$ we know by \ref{subsec:basic-properties-of-extension}
that $M_{i}\underset{I}{\overset{f_{i}}{\longrightarrow}}M_{i+1}$. 

if $i>1$, for every $0<j<i$ the fact that $y_{j}$ calculates $x_{j}$
is also absolute in sub-models (by \ref{subsec:absoluteness-of-Turing})
and so within $M_{i+1}$ every $f_{i}\left(f_{i-1}\left(\cdots\left(f_{j}\left(y_{j}\right)\right)\right)\right)$
calculates $f_{i}\left(f_{i-1}\left(\cdots\left(f_{j}\left(x_{j}\right)\right)\right)\right)$
for every $0<j<i$. As $y_{j}<\mathcal{L}\left(T\right)$ (and $\mathcal{L}\left(T\right)$
is absolute), in $M_{j+1}$ we know that 
\[
f_{i}\left(f_{i-1}\left(\cdots\left(f_{j}\left(y_{j}\right)\right)\right)\right)<f_{i}\left(f_{i-1}\left(\cdots\left(f_{0}\left(\mathcal{L}\left(T\right)\right)\right)\right)\right)
\]

Now fix $i=\mathcal{L}\left(T\right)+2$. In $M_{i+1}$ for such $i$
exist 
\[
y_{0},y_{1},\ldots,y_{i}
\]
 each of them smaller than $\mathcal{L}\left(T\right)$ and each $y_{j}$calculates
$x_{j}$. As all 
\[
x_{0},x_{1},\ldots,x_{i}
\]
were different we get that 
\[
y_{0},y_{1},\ldots,y_{i}
\]
 are all different. So the model $M_{i+1}$has $i=\mathcal{L}\left(T\right)+2$
different numbers all smaller than $\mathcal{L}\left(T\right)$, which
is a contradiction. Therefore, our assumption that the claim $con\left(T\right)$
could be proven from $T$ was wrong $\blacksquare$

\end{document}